\documentclass{article}
\usepackage[utf8]{inputenc}
\usepackage{amsmath,amsthm,amssymb}
\usepackage{authblk}
\usepackage{graphicx}
\usepackage{enumerate}
\usepackage{hyperref}
\usepackage{fullpage}

\hypersetup{
    colorlinks=true,
    linkcolor=blue,
    pdftitle={A sharper Ramsey theorem for constrained drawings},
    pdfauthor={Pavel Paták}
    }

\begin{document}

\title{A sharper Ramsey theorem for constrained drawings}

\author[1,2]{Pavel Paták}

\affil[1]{Department of Applied Mathematics, Faculty of Mathematics and Physics, Charles University, Czech Republic}

\affil[2]{Department of Applied Mathematics, Faculty of Information Technology, Czech Technical University in Prague,  Czech Republic}




\newtheorem{conjecture}{Conjecture}
\newtheorem{theorem}{Theorem}
\newtheorem{proposition}{Proposition}
\newtheorem{lemma}{Lemma}
\newtheorem{definition}{Definition}

\maketitle
\abstract{Given a graph $G$ and a collection $\mathcal C$ of subsets of $\mathbb{R}^d$ indexed by the subsets of vertices of $G$, a constrained drawing of $G$ is a drawing, where each edge is drawn inside some set from $\mathcal C$, in such a way that non-adjacent edges are drawn in sets with
disjoint indices. In this paper we prove a Ramsey type result for such drawings.
Furthermore we show how the result can be used to obtain
Helly type theorems.

More precisely, we prove the following.
For each $n$ and $b$, there is $N=O(b^{2n-3})$ with the following properties:
If $G$ is a drawing of a graph on $N$ vertices and $\mathcal C$ is a collection
of sets of $\mathbb{R}^d$ such that each $(b+1)$-tuple $T$ of vertices lies in a set
indexed by $T$ and contains at least one edge in $T$,
then in $G$, we can find a constrained copy of the complete graph $K_n$.

As a direct consequence we obtain the following Helly type result:
For each $d$, there is a polynomial $h(b)$ of degree at most $2d+3$ such that the following holds.
For every family $\mathcal F$ of sets in $\mathbb{R}^d$, its Helly number is at most $h(b)$,
provided that the intersection of any non-empty subfamily has at most $b$
path-connected components, and trivial homology groups $H_1$, $H_2$, .... $H_{\lceil d/2\rceil-1}$.
This dramatically improves the original theorem by Matoušek which
had stronger assumption and a tower-like bound on $h(b)$.
Under the same assumptions, our technique can also be used to bound
Radon numbers.}





\section{Introduction}
The introduction of Ramsey theorems for constrained drawings has spurred the development of topological Helly type results in the recent years~\cite{m-httucs-97,hb17,patakova2020,holmsen2021}. 
However, the bounds obtained in these theorems are too large for practical applications.
In the present paper we provide the first step to make the bounds practical: we focus on constrained graph drawings and provide a polynomial upper bound for the corresponding Ramsey problem. Consequently, we obtain a polynomial bound for some of the Helly type results.

In order to put the definitions and the main theorem into the proper context, we first explain the applications.
\begin{definition}
Let $X$ be a set. 
A \emph{closure operator on $X$} is any map $\operatorname{cl}\colon 2^X\to 2^X$ that
satisfies the following three conditions for all $A,B\subseteq X$.
\begin{description}
 \item[Extensivity] $A\subseteq \operatorname{cl}A$ 
 \item[Idempotence] $\operatorname{cl}\left(\operatorname{cl}A\right)=\operatorname{cl}A$
 \item[Monotonicity] $A\subseteq B\Rightarrow \operatorname{cl}A\subseteq\operatorname{cl}B$.
\end{description}
\end{definition}
The reader may view closure operators as a generalization of convex and affine hulls.
However, there are also exotic closure operators for which $\operatorname{cl}\emptyset\neq\emptyset$. 

A closure operator may be identified with an intersection closed family of sets:
If $\operatorname{cl}$ is a closure operator on $X$, then the sets $A$ with $\operatorname{cl}A=A$ form an intersection closed family. Conversely, if $\mathcal F$ is any family of subsets of a given set $X$, then we may define the closure operator $\operatorname{cl}_{\mathcal F}$ via
\[
\operatorname{cl}_{\mathcal F}(S):=\bigcap \{F\in\mathcal F\mid S\subseteq F\}. 
\]
If there is no set $F\in\mathcal F$, the closure $\operatorname{cl}_\mathcal F(S)$ is defined to be $X$. Conversely, if $\operatorname{cl}$ is any closure operator on $X$, \[\mathcal F_{\operatorname{cl}}:=\{S\subseteq X\mid \operatorname{cl}S=S\}\] is a intersection closed family of sets with $\operatorname{cl}_{\mathcal F_{\operatorname{cl}}}=\operatorname{cl}$. This allows us to pass from closure operators to intersection closed set families and back.

\begin{definition}
For a closure operator $\operatorname{cl}$ on $X$, its \emph{Radon number $r\left(\operatorname{cl}\right)$} is defined as the smallest integer $r$ such that any set $S\subseteq X$ of size at least $r$ 
contains two disjoint subsets $S_1,S_2\subseteq S$ with $\operatorname{cl} S_1\cap \operatorname{cl} S_2\neq\emptyset$.  If no such integer exists, we set $r(\operatorname{cl})=\infty$.
For a family $\mathcal F$ of sets,
we define $r(\mathcal F):=r(\operatorname{cl}_\mathcal F)$.
\end{definition}
Any upper bound on Radon number is a strong result, as it implies a Helly-type result~\cite{Levi1951}, a version of Tverberg's theorem~\cite{jamison1981} and plethora of other convexity results for $\mathcal F$, including variants of fractional Helly theorem, colorful Helly theorem, existence of weak $\varepsilon$-nets and $(p,q)$-theorem~\cite{holmsen2021}.

\begin{figure}
  \begin{center}
    \includegraphics[page=5]{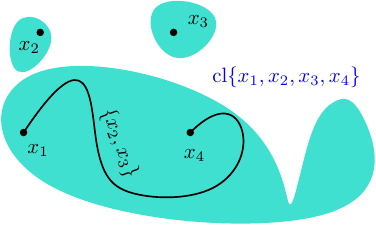}
    \caption{An edge labeled $\{x_2,x_3\}$ connecting $x_1$ to $x_2$ inside $\operatorname{cl}\{x_1,x_2,x_3,x_4\}$.}
    \label{fig:connection}
  \end{center}
 \end{figure}

The papers by Matoušek~\cite{m-httucs-97}, its generalization by Goaoc, Paták, Patáková, Tancer and Wagner~\cite{hb17} and the subsequent work of Patáková~\cite{patakova2020}
introduced a powerful technique that uses constrained drawings to bound Radon numbers.
Let us explain the main idea of the technique on a simple example.
Let $\operatorname{cl}$ be a closure operator in $\mathbb{R}^2$ such that for every $F\subseteq \mathbb{R}^2$, $\operatorname{cl} F$ has at most three path-connected components.
Given a set $S$ of $N$ points in $\mathbb{R}^2$, we know that if we take a set $F\subseteq S$ of four of these points, two points from $F$ can be connected with a path that lies within $\operatorname{cl}F$.
In such a case, we draw this edge $e$ and assign it the label $F\setminus e$, see~Figure \ref{fig:connection}.
Suppose now that for each set $S$ on $N$ points we find a copy of $K_5$,
where disjoint edges have disjoint labels.
Then, by the Hanani-Tutte theorem~\cite{Chojnacki1934,Tutte70}, there are two non-adjacent edges of this $K_5$ whose drawings intersect. Without loss of generality, let these edges be $v_1v_2$ with label $\{v_3,v_4\}$
and $v_5v_6$ with label $\{v_7,v_8\}$ and let the intersection point be $p$.
But then $p\in \operatorname{cl}\{v_1,v_2,v_3,v_4\}\cap\operatorname{cl}\{v_5,v_6,v_7,v_8\}$.
Consequently, $r(\operatorname{cl})\leq N$. 

Therefore, in order to establish good bounds on the Radon number for sets in $\mathbb{R}^2$,
we need to establish a value of $N$ that will guarantee the desired copy of $K_5$ (or $K_{3,3}$). For sets in higher dimensions, we shall establish a value  of $N$ that guarantees a copy of $K_n$.

In our main theorem (Theorem~\ref{thm:bounds}), we construct polynomials $p_n(b)=O(b^{2n-3})$, such that if $\operatorname{cl} Y$ has at most $(b+1)$ path-connected components for any $Y$, then in order to find a copy of $K_n$ with disjoint edge labels, it suffices to take $N\geq p_n(b)$. This greatly improves upon the previous tower function bound by Matoušek~\cite[Theorem 2]{m-httucs-97}.

\paragraph*{Notation}
If $P$ is a set, the symbol $2^P$ stands for the family of all its subsets. The symbol $\binom{P}{k}$ denotes the set of all $k$-element subsets of $P$. $\mathbb{Z}_2$ denotes the two-element field. 
An (abstract) simplicial complex $K$ is a family of finite sets, that is closed with respect to taking subsets. The elements of $K$ are called simplices. If a simplex $\sigma\in K$ has $n$-elements, we say that its dimension is $n-1$. 
The $n$-simplex $\Delta_n$ is the abstract simplicial complex formed by all subsets of $\{1,2,\ldots, n, n+1\}$. If $K$ is an abstract simplicial complex, its $k$-skeleton  $K^{(k)}$ is the simplicial complex consisting of all $\leq k$-dimensional simplices of $K$.
If we write that $\Phi\colon K\to 2^P$ is a map, we mean that $\Phi$ assigns a subset of $P$ to every simplex in $K$. If $K\subseteq 2^{\{1,2,\ldots, n\}}$ is an abstract simplicial complex, its geometric realization $|K|$ is a topological space obtained as follows:
We consider the standard basis $e_1,\ldots, e_n$ of $\mathbb{R}^n$. For a simplex $\sigma\in K$, its geometric realization $|\sigma|$ equals $\operatorname{conv} \{e_i\mid  i\in\sigma\}$, and $|K|:=\bigcup_{\sigma\in K} |\sigma|$.

The symbol $\sqcup$ stands for disjoint union. A graph $G$ is a pair $(V,E)$, where $V$ and $E$ are disjoint sets and $E\subseteq  \binom{V}{2}$, the elements of $V$ are called the vertices of $G$ and the elements of $E$ are its edges. If $G$ is a graph, the symbols $V(G)$ and $E(G)$ stand for the set of vertices and the set of edges of $G$, respectively.
$K_n$ is the complete graph on $n$ vertices, $K_{m,n}$ is the complete bipartite graph with parts of sizes $m$ and $n$. 
If $G$ and $H$ are two graphs, $G+H$ denotes their graph join, that is a graph obtained from their disjoint union by adding all the edges $\{g,h\}, g\in V(G), h\in V(H)$. Note that the operation of graph join is commutative and associative. If $t\in\mathbb{N}\cup\{0\}$ and $G$ is a graph, $t\cdot G$ stands for the disjoint union of $t$ copies of $G$.

The symbol $H_i(X;\mathbb{Z}_2)$ denotes the homology group of $X$ with $\mathbb{Z}_2$ coefficients,
the symbol $\widetilde H_i(X;\mathbb{Z}_2)$ denotes the reduced homology group of $X$ with $\mathbb{Z}_2$ coefficients. By definition, if $X$ is non-empty, then $H_0(X;\mathbb{Z}_2)$ is isomorphic to $\mathbb{Z}_2^{\# \text{of path-components of }X}$ and $\widetilde{H_0}(X;\mathbb{Z}_2)$ is isomorphic to $\mathbb{Z}_2^{\# \text{of path-components of }X - 1}$ and $H_i(X;\mathbb{Z}_2)=\widetilde{H}_i(X;\mathbb{Z}_2)$ for all $i\geq 1$.
In particular if $X$ is contractible, then $\widetilde{H}_i(X;\mathbb{Z}_2)=0$ for all $i$.
The symbol $\pi_i(X)$ denotes the $i$th homotopy group. In particular, if $B$ is the $(i+1)$-dimensional ball, then $\pi_i(X)=0$ if every map $f\colon\partial B\to X$ can be extended to a map $f\colon B\to X$.

The paper is organized as follows. In Section~\ref{sec:main} we introduce the necessary definitions and state our main results. The proofs of these results are then showed in Sections~\ref{sec:combinatorics} and \ref{sec:radon}. In Section~\ref{sec:combinatorics} we focus on the combinatorial part, prove the main theorem and discuss some improvements and optimality of the results. The proofs in this sectional are purely combinatorial.
In Section~\ref{sec:radon}, we prove the Helly type theorems from Section~\ref{sec:main} by combining our main theorem with the technique of constrained drawings. The proofs in this section rely on some arguments from algebraic topology. 
In Section~\ref{sec:open}, we list several open problems and conjectures.

\section{Main results}
\label{sec:main}
\begin{figure}
 \begin{center}
 \includegraphics[page=2]{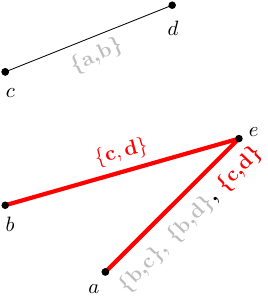}
 \caption{A graph $G$ with edge labels and a highlighted constrained copy of $K_{1,2}$}\label{fig:constrained}
 \end{center}
\end{figure}

\begin{definition}
 Let $H$ be a graph with vertex set $V(H)$,
 where each edge $e\in E(H)$ is assigned a set of labels $L_e\subseteq 2^{V(H)\setminus e}$.
 We say that $H$ \emph{contains a constrained copy of a graph $G$} if the following holds.
 \begin{enumerate}
  \item There is a subgraph $G_0$ of $H$ that is isomorphic to $G$.
  \item For this subgraph $G_0$, there is an assignment of edge labels $l_0\colon E(G_0)\to 2^{V(G)}$ such that 
  \begin{enumerate}
  \item for each $e\in E(G_0)$, $l_0(e)\in L_e$,
  \item for each $e\in E(G_0)$, $l_0(e)$ is disjoint from $V(G_0)$
  \item if $e,f\in E(G_0)$ are disjoint,
  then $l_0(e)$ and $l_0(f)$ are disjoint.
  \end{enumerate}
 \end{enumerate}
\end{definition}
See Figure~\ref{fig:constrained} for an illustration of the notion.
\begin{definition}
 Let $b\geq 1$ be an integer.
 We say that a graph $H$ is \emph{$b$-iatlon}, if each edge $e\in E(H)$ of $H$ is assigned a set of labels $L_e\subseteq 2^{V(H)\setminus e}$ and 
 if for each subset $S\subseteq V(H)$
 with $|S|=b+1$, there exist $u,v\in S$
 for which $\{u,v\}$ is an edge of $H$
 and the set $S\setminus\{u,v\}$ lies in $L_{\{u,v\}}$.
\end{definition} 
The term $b$-iatlon is an abbreviation from \textbf{i}n \textbf{a}ll \textbf{t}uples \textbf{l}abeled \textbf{on}e. Note that the graph $G$ in Figure~\ref{fig:constrained} is $3$-iatlon. Figure~\ref{fig:connection} provides the motivation why to study $b$-iatlon graphs and constrained copies. Note that every graph with $\leq b$ vertices is $b$-iatlon.

We are now finally in a position to state our main theorem that every sufficiently large $b$-iatlon graph contains a constrained copy of $K_n$.
\begin{theorem}\label{thm:bounds}
 Let $k,m,n,b$ be non-negative integers with $1\leq m\leq n$, $2\leq n$ and $1\leq b$. For simplification, for each non-negative integer $l$, we set
 $S_l:=\sum_{j=0}^{l} \binom{b+1}{2}^j$. 
 Let $H$ be a $b$-iatlon graph with $N$ vertices. Then
 \begin{enumerate}
  \item If $N\geq b\cdot S_{n-2} + 1$, then $H$ contains a constrained copy of $K_n$.  \label{it:five}
  \item If $N\geq b\cdot S_{n-2} + 1 + (k-1)(bn-b+1)$, then $H$ contains a constrained copy of $k\cdot K_n$.\label{it:kfive}
  \item If $N\geq \binom{b+1}{2}^m\cdot (n-1) + b\cdot S_{m-1}+1$, then $H$ contans a constrained copy of $K_{m,n}$.\label{it:bipartite}
  \item 
  If $N\geq \binom{b+1}{2}^m(n-1)+b\cdot S_{m-1}+1 + (k-1)(bm+n)$, then $H$ contains a constrained copy of $k\cdot K_{m,n}$.
  \label{it:kbipartite}
 \end{enumerate}
\end{theorem}
We note that if $b>1$, $S_l=\left(\binom{b+1}{2}^{l+1}-1\right)/\left(\binom{b+1}{2}-1\right)$. Unfortunately, this does not work if $b=1$.

Before we prove a stronger version of the theorem in Section~\ref{sec:combinatorics}, let us look at several of its corollaries.

If we apply the result to constrained drawings in Patáková's proof of topological Radon's theorem~\cite{patakova2020}, we improve her upper bound on $r(\operatorname{cl})$.
Under some additional assumption we reduce the tower function bound to a polynomial one. 
\begin{theorem}\label{thm:radon}
 Let $\operatorname{cl}$ be a closure operator on $\mathbb{R}^d$. Let $b\in\mathbb{N}$ satisfy
 \begin{enumerate}
  \item  for every $X\subseteq\mathbb{R}^d$ with $|X|=b+1$, $\operatorname{cl} X$ has at most $b$ path-connected components; and
  \item for any integer $i$ satisfying $1\leq i\leq \lceil d/2\rceil-1$ and any $Y\subseteq \mathbb{R}^d$ with \[b+1+i\leq |Y|\leq i+2+(b-1)\binom{i+2}{2}\text{, one has }\widetilde{H}_{i}(\operatorname{cl} Y;\mathbb{Z}_2)=0.\] 
 \end{enumerate}
 Then 
 for $S_{d+1}=\sum_{j=0}^{d+1}\binom{b+1}{2}^j$
 \[r(\operatorname{cl})\leq b\cdot S_{d+1} + 1=O(b^{2d+3}).\]
\end{theorem}
Moreover, if $d=2$, the bound can be improved to $r(\operatorname{cl})=O(b^6)$,
see~Theorem~\ref{thm:surfaces}.

As an immediate consequence, we obtain polynomial bounds on so called Helly numbers.
\begin{definition}
 Let $\mathcal F$ be a family of sets.
 \emph{Helly number $h(\mathcal F)$} is the smallest number $h\in\mathbb{N}\cup\{\infty\}$ such that all finite subfamilies $\mathcal G\subseteq\mathcal F$ with $\bigcap\mathcal G=\emptyset$ contain $h$ (not necessarily distinct) sets $G_1,G_2,\ldots, G_h\in\mathcal G$ with $G_1\cap\cdots \cap G_h=\emptyset$.
\end{definition}
For many algoritmic applications, Helly number measures how much a given family $\mathcal F$ of sets behaves like a family of convex sets in $\mathbb{R}^d$. This is also a reason why it has been thoroughly studied~\cite{danzer1963,eckhoff93,wenger97,tancer2013,amenta17} and establishing good upper bounds under various conditions is a fundamental question.

Since by Levi~\cite{Levi1951}, $h(\mathcal F)+1\leq r(\mathcal F)$, 
Theorem~\ref{thm:radon} immediately implies the following bound. 
\begin{theorem}\label{thm:helly}
 Let $\mathcal F$ be a family of sets in $\mathbb{R}^d$. Let $b\in\mathbb{N}$ satisfy
 \begin{enumerate}
  \item for every finite subfamily $\mathcal H\subseteq\mathcal F$, the intersection $\bigcap\mathcal H$ has at most $b$ path-connected components; and
  \item for every finite $\mathcal H\subseteq\mathcal F$ and every $i=1,2,\ldots,\lceil d/2\rceil-1$, 
  $\widetilde{H}_i\left(\bigcap\mathcal H;\mathbb{Z}_2\right)=0$.
 \end{enumerate}

 Then for $S_{d+1}=\sum_{j=0}^{d+1}\binom{b+1}{2}^j$
 \[h(\mathcal F)\leq  b\cdot S_{d+1} = O(b^{2d+3}).\]
\end{theorem}
The first result of this type was proven by Matoušek~\cite{m-httucs-97}. However, he used stronger assumption $\pi_i(\bigcap \mathcal H)=0$ and his proof method led to an enormous bound $h(\mathcal F)\leq R_{b+1}\left(\binom{2d+3}{2}(b-1) + 2d+3; \binom{b+1}{2}\right)$, where $R_k(n;c)$ denotes the hypergraph Ramsey number. 
We recall that $R_k(n;c)$ is the smallest $N$ such that every coloring of $k$-element subsets of $N$ with $c$ colors will contain a monochromatic set of size $n$. For practical purposes such a bound is too huge, as for $k\geq 4$ and $c\geq 2$, $R_k(n,c)\geq R_k(n,2)\geq \operatorname{twr}_{k-1}(Cn^2)$, where $\operatorname{twr}_k(x)$ denotes the tower function $\operatorname{twr}_1(x):=x$, $\operatorname{twr}_{l+1}(x):= 2^{\operatorname{twr}_l(x)}$, and $C$ is some suitable constant~\cite{graham1990}.
Theorem~\ref{thm:helly} is hence a significant improvement of Matoušek's result.

We remark that later on Matoušek's result has been generalized by Goaoc, Paták, Patáková, Tancer and Wagner~\cite{hb17} into the following form:
There exists a number $h'(b,d)$ such that the following holds.
If $\mathcal F$ is a family of sets in $\mathbb{R}^d$ such that for all $i=0,1,\ldots, \lceil d/2\rceil-1$, $\dim \widetilde H_i(\mathcal F;\mathbb{Z}_2)\leq b$, then $h(\mathcal F)\leq h'(b,d)$. The proof uses the hypergraph Ramsey theorem iteratively, and so the obtained bound $h'(b,d)$ is even larger than in the case of Matoušek's theorem~\cite{m-httucs-97}. For this generalized theorem, our proof method can be used to obtain better bounds for the first induction step. This effectively peels off several levels of the resulting tower functions, but still leaves the bounds impractically large.
When one tries to generalize our approach to the other induction steps in order to obtain good bounds, one encounters substantial obstacles. Overcoming these barriers seems to require several novel ideas.

By~Levi~\cite{Levi1951}, $h(\mathcal F)+1\leq r(\mathcal F)$, hence from now on we focus on bounding Radon numbers only.
In the next theorem, we bound Radon number of closure operators on surfaces.
\begin{theorem}\label{thm:surfaces}
 Let $\operatorname{cl}$ be an arbitrary closure 
 operator on a closed surface $M$.
 Let $g$ be the genus of $M$, if $M$ is orientable, and the non-orientable genus of $M$, if $M$ is non-orientable.
 
 Let $b$ be a positive integer such that
 for every $X\subseteq M$ of cardinality $b+1$, $\operatorname{cl} X$ has at most $b$ path-connected components. Set $x:=\binom{b+1}{2}$. Then
 \[
  r(\operatorname{cl})\leq 2x^3 + bx^2 + bx + b + 3g(b+1)=2x^3 + bx^2 + bx + O(b).
 \]
 where the constant in $O$-notation depends on $g$.
\end{theorem}
More generally, the technique
shows that topological conditions
guarantee a bound on Radon number
for closed compact PL manifolds.
\begin{theorem}\label{thm:manifolds}
 Let $M$ be a closed compact PL manifold of dimension $2d$. Let $\beta_d(M;\mathbb{Z}_2)$ be the middle Betti number of $M$.
 Let $\operatorname{cl}$ be an arbitrary closure operator on $M$ and let $b$ be a positive integer.
 Let $t=2d+1 + (d+1)\beta_d(M;\mathbb{Z}_2)$.
 
 Suppose that
 \begin{enumerate}
  \item For each $X\subseteq M$ with $|X|=b+1$, $\operatorname{cl} X$ has at most $b$ path-connected components and 
  \item for each $i=1,2,\ldots, \lceil d/2\rceil-1$ and each finite $X\subseteq M$ with \[b+1+i\leq |X|\leq i+2+(b-1)\binom{i+2}{2},\quad \pi_i(\operatorname{cl} X)=0,\]
 \end{enumerate}
 then $r(\operatorname{cl})\leq b\cdot S_{t+1}+1$, where $S_{t+1}=\sum_{j=0}^{t+1} \binom{b+1}{2}^j$.
 \end{theorem}
 
 The proofs of Theorems~\ref{thm:radon}, \ref{thm:surfaces} and \ref{thm:manifolds} can be found in Section~\ref{sec:radon}.
 
 \section{Combinatorial part}\label{sec:combinatorics}
 In this section we prove somewhat stronger version of Theorem~\ref{thm:bounds} and introduce all the definitions and tools that are needed to do so.
 The most versatile of these tools is going to be Proposition~\ref{prop:sRamsey}, which allows us to find the desired constrained copy of a graph $G$ by using induction on $|V(G)|$. The argument relies on finding constrained copies of stars $K_{1,n}$. The optimality of the results for stars is discussed in Subsection~\ref{sec:optimality}.
 This section is completely combinatorial and does not rely on any topological notions.
 
 \begin{definition}
 Let $H$ be a graph with vertex set $V(H)$, where each edge $e\in E(H)$ is assigned a set of labels $L_e\subseteq 2^{V(H)\setminus e}$.
 Let $X\subseteq V(H)$ be a set.
 We define the \emph{restriction $H[X]$ of $H$ to $X$} as follows: as a graph, it is the subgraph of $H$ induced by $X$; as for the labels, for each edge $e\in E(H)$, we only keep those labels $l_e\in L_e$ for which $l_e\subseteq X$.
\end{definition}
Cleary, if $H$ is a $b$-iatlon graph and $X\subseteq V(H)$ a subset of its vertices, the restriction $H[X]$ is also a
$b$-iatlon graph. 

Our aim is to show that for every graph $G$ there is some $n$ such that every $b$-iatlon graph with at least $n$ vertices contains a constrained copy of $G$. In order to obtain better bounds, we shall also keep track of the number of distinct edge labels in the desired constrained copy. Therefore, we define the following two  functions.
\begin{definition}
 Let $G$ be a graph and $s,b\in\mathbb{N}$ be positive integers. Then we set
 \begin{align*}
 r_s(G,b)&:=\min \{n\mid \text{any $b$-iatlon graph on $n$ vertices contains a constrained copy of $G$ with $\leq s$ distinct edge labels}\}\\
 p(G,b)&:=\min \{n\mid \text{any $b$-iatlon graph on $n$ vertices contains a constrained copy of $G$}\}
  \end{align*}
  If the minimal value does not exist, we define the functions to be $+\infty$.
\end{definition}
It is easy to observe that these functions have the following properties:
\begin{enumerate}
 \item The function $r_s(G,b)$ is non-increasing in $s$: if $s<t$,
then $r_s(G,b)\geq r_t(G,b)$.
 \item If $s\geq |E(G)|$, then $r_s(G,b)=r_{|E(G)|}(G,b)$
 \item $r_{|E(G)|}(G,b)=p(G,b)$.
\end{enumerate}
Our goal is to show that for every graph $G$, the function $p(G,b)$ is bounded from above by some polynomial in $b$. For that we provide a polynomial upper bound for $r_{|V(G)|-1}(G,b)\geq r_{|E(G)|}(G,b)=p(G,b)$.

\begin{proposition}\label{prop:sRamsey}
 For all finite graphs $G$ and $L$ and all positive integers $b,s,t$ the following holds.
 \begin{enumerate}
  \item $r_0(K_1,b)=1$. \label{it:vertex}
  \item $r_1(K_2,b)=b+1$. \label{it:edge}
  \item $r_1(G,1)=|V(G)|$. \label{it:noholes}
  \item $r_{s+t}(G\sqcup L,b)\leq r_s(G,b)+r_t(L,b)$.\label{it:dunion}
  \item $r_1(K_{1,n},b)\leq \binom{b+1}{2}(n-1)+b+1$.\label{it:star}
  \item $r_{s+t}(K_1+G,b)\leq
  r_s(K_{1,r_t(G,b)},b)$.\\
  In particular, 
  $r_{1+t}(K_1+G,b)\leq 
  \binom{b+1}{2}(r_t(G,b)-1)+b+1$.\label{it:gcone}
  \item If $G'$ is a subgraph of $G$, then 
  $r_s(G',b)\leq r_s(G,b)$ for all positive integers $b,s$.\label{it:subgraph}
  \item If $G$ has at least two vertices, then $r_{|V(G)|-1}(G,b)=O\left(b^{2\cdot |V(G)|-3}\right)$.\label{it:degree}
  \item If $|V(G)|+(b-1)\cdot s + r_t(L,b)\geq r_s(G,b)$,
  then $r_{s+t}(G\sqcup L,b) \leq |V(G)|+(b-1)\cdot s + r_t(L,b)$.
  \label{it:bdunion}
 \end{enumerate}
\end{proposition}

\begin{proof}
Note that all the inequalities of the theorem are automatically satisfied, if the right side equals $+\infty$. Therefore, during the whole proof we assume that this is not the case.
 \begin{description}
  \item[(\ref{it:vertex})] Since $K_1$ does not contain any edge, any copy of $K_1$ is automatically constrained with no edge labels.
  \item[(\ref{it:edge})]
  If $H$ is a $b$-iatlon graph with $b+1$ vertices, it contains an edge with some label and this edge is a constrained copy of $K_2$ with a single label. The graph on $b$ vertices without any edges is $b$-iatlon and does not contain $K_2$.
  \item[(\ref{it:noholes})] If $H$ is a $1$-iatlon graph on $n=|V(G)|$ vertices, it is a complete graph and for each of its edges, $\emptyset$ is one of its labels.
  Therefore, $H$ contains a constrained copy of $G$ with all edge labels $\emptyset$.
  \item[(\ref{it:dunion})] 
  Let $H$ be a $b$-iatlon graph on $r_s(G,b)+r_t(L,b)$ vertices. We split the vertices of $H$ into two parts $P_G$ and $P_L$ of sizes $r_s(G,b)$ and $r_t(L,b)$, respectively.
  By the assumption, $H[P_G]$ contains a constrained copy $C_0$ of $G$ with $\leq s$ distinct edge labels and $H[P_L]$ contains a constrained copy $C_1$ of $L$ with $\leq t$ distinct edge labels.
  Then $C_0\sqcup C_1$ is a constrained copy of $G\sqcup L$ in $H$ with at most $s+t$ distinct edge labels.
  \item[(\ref{it:star})]  
  Let $H$ be a $b$-iatlon graph with $N$ vertices. Let us double count the number $T$ of triples $(v,e,l)$, where $v$ is a vertex of $H$, $e$ is an edge incident to $v$ and $l\in L_e$ is a label of $e$ with $|l|=b-1$.
  Since $H$ is $b$-iatlon, each $(b+1)$-tuple of vertices contains at least one labeled edge with $|l|=b-1$, and this edge has two endpoints.
  Hence $2\binom{N}{b+1}\leq T$. 
  Let $m$ be the minimal non-negative integer such that for each $v\in V(H)$, at most $m$ edges incident to $v$ have the same label.
  Then $T\leq N\cdot \binom{N-1}{b-1}m$, since we can first choose the vertex $v$, then the possible label $l$ (which does not contain $v$), and then there are at most $m$ edges with label $l$ that are incident to $v$.
  Therefore,
   $2\binom{N}{b+1} \leq N\cdot \binom{N-1}{b-1}m$, which simplifies to
   \begin{equation}\label{eq:degreeBound}
    2\frac{N-b}{(b+1)b} \leq m.
   \end{equation}
   If we plug in the value $N=\binom{b+1}{2}(n-1)+b+1$, we  get
   $\frac{\binom{b+1}{2}(n-1)+1}{\binom{b+1}{2}}\leq m.$
   Since $m$ is an integer, $m\geq n$.
   But this means that there is a vertex $v$ with $n$ adjacent edges which all have the the same label $l$.
    By definition of the labels, such $l$ does not contain $v$, nor any endpoint of the $n$ edges. Therefore, these $n$ edges together with their endpoints and the label $l$ form a constrained copy of $K_{1,n}$ in $H$ and all edges of this constrained copy use the same label.
  \item[(\ref{it:gcone})] 
  If $H$ is a $b$-iatlon graph on $r_s(K_{1,r_t(G,b)},b)$ vertices, then we can find a constrained copy $C_0$ of $K_{1,r_t(G,b)}$ inside $H$ with $\leq s$ distinct edge labels. Let $v_0$ be the root of this star and let $X$ be the set of its leaves.
  By the definition of $r_t(G,b)$, there is a constrained copy $C_1$ of $G$ in $H[X]$ which uses $\leq t$ labels,
  see~Figure~\ref{fig:induction}.
  \begin{figure}
   \begin{center}
   \includegraphics[page=3]{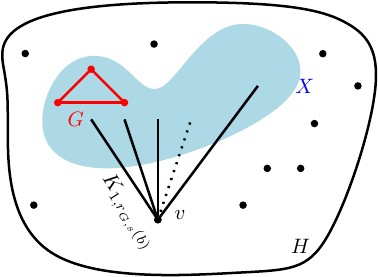} 
   \caption{Finding $K_1+G$}\label{fig:induction}
   \end{center}
  \end{figure}
  Then, since the edge labels of $C_0$ are disjoint from $X$ and $v_0$, $C_0\cup C_1$ is a constrained copy of $G\cup K_{1,r_s(G,b)}$ which uses $\leq (s+t)$ distinct edge labels.
  If we now disregard the leaves of $C_0$ that do not belong to $C_1$, and consider the restriction of $C_0\cup C_1$ to $\{v_0\}\cup C_1$, we obtain a constrained copy of $K_1+G$ with $\leq (s+t)$ distinct edge labels. 
  
  Plugging in
  $r_1(K_{1,n},b)\leq \binom{b+1}{2}(n-1)+b+1$ from \ref{it:star},
  yields $r_{1+t}(K_1+G,b)\leq 
  \binom{b+1}{2}(r_t(G,b)-1)+b+1$, as desired.
  
  \item[(\ref{it:subgraph})] This follows from the fact that a subgraph of a constrained copy is again constrained.
  \item[(\ref{it:degree})] Due to (\ref{it:subgraph}), it suffices to prove the statement for complete graphs $K_n$.
  We proceed by induction. For $n=2$, the statement is true by (\ref{it:edge}).
  For $n\geq 2$, (\ref{it:gcone}) yields
  $r_n(K_{n+1},b)=r_{n}(K_1+K_n,b)\leq \binom{b+1}{2}(r_{n-1}(K_n,b)-1)+b+1=O\left(b^{2+2n-3}\right) = O\left(b^{2n-1}\right)$.
  \item[(\ref{it:bdunion}] 
  Let $H$ be a $b$-iatlon graph on $|V(G)|+(b-1)s + r_t(L,b)\geq r_s(G,b)$ vertices. By definition of $r_s(G,b)$, $H$ contains a constrained copy $G_0$ of $G$ with $\leq s$ distinct edge labels. Let $X$ be the union of $V(G)$ with all the labels of $G_0$. Then, since each label has size $(b-1)$, $|X|\leq |V(G)|+(b-1)s$.
  Therefore, the restriction of $H$ to $V(H)\setminus X$ is a $b$-iatlon graph with at least $r_t(L,b)$ vertices, and contains a constrained copy $L_0$ of $L$ with $\leq t$ distinct edge labels.
  Then $G_0\sqcup L_0$ is a constrained copy of $G\sqcup L$ with $\leq (t+s)$ distinct edge labels, as desired.\qedhere
 \end{description}
\end{proof}
Further improvements of Proposition \ref{prop:sRamsey} are discussed in subsection \ref{sec:optimality}.

We are now ready to prove the following stronger version of Theorem~\ref{thm:bounds}, which also indicates how many distinct edge labels we use. 
\begin{theorem}\label{thm:boundsColors}
 Let $k,m,n,b$ be non-negative integers with $1\leq m\leq n$, $2\leq n$ and $1\leq b$. For each non-negative integer $l$, we set
 $S_l:=\sum_{j=0}^{l} \binom{b+1}{2}^j$. 
 Then
 \begin{enumerate}
  \item $r_{n-1}(K_n,b)\leq bS_{n-2}+1$, \label{it:Bfive}
  \item $r_{k(n-1)}(k\cdot K_n,b)\leq S_{n-2} + 1 + (k-1)(bn-b+1)$,\label{it:Bkfive}
  \item $r_m(K_{m,n},b)\leq \binom{b+1}{2}^m(n-1)+b\cdot S_{m-1}+1$,\label{it:Bbipartite}
  \item $r_{km}(k\cdot K_{m,n}\leq \binom{b+1}{2}^m(n-1)+b\cdot S_{m-1}+1 + (k-1)(bm+n)$.
  \label{it:Bkbipartite}
  \end{enumerate}
\end{theorem}
Since for all $s$, $p(G,b)\leq r_s(G,b)$, the result immediately yields Theorem~\ref{thm:bounds}.

\begin{proof}
 \begin{description}
  \item[(\ref{it:Bfive})] 
  We proceed by induction on $n$. By Proposition~\ref{prop:sRamsey}(\ref{it:edge}), $b+1 = b\cdot S_0+1$ vertices of $H$ suffice to force a constrained copy of $K_2$, which has necessarily a single edge label. 
  So suppose that the statement holds for $n$. Then by Proposition~\ref{prop:sRamsey}(\ref{it:gcone}), a constrained copy of $K_{n+1}=K_1+K_n$ with $\leq n$ distinct edge labels is guaranteed as soon as
  $H$ has at least $\binom{b+1}{2}\cdot (b\cdot S_{n-2}+1-1) + b + 1= b\cdot S_{n-1}+1$ vertices.
  \item[(\ref{it:Bkfive})] We proceed by induction on $k$, the case $k=1$ being given by the previous point.
  Since \[b\cdot S_{n-2}+1+(k-1)(bn-b+1)\geq b\cdot S_{n-2}+1\geq r_{n-1}(K_n,b),
  \]
  Proposition~\ref{prop:sRamsey}(\ref{it:bdunion}) and the induction assumption yield the desired inequality
  \begin{align*}r_{k(n-1)}(k\cdot K_n,b)
  &\leq n + (b-1)(n-1) + r_{(k-1)(n-1)}((k-1)\cdot K_n,b)\\&\leq 
  (bn-b+1) + b\cdot S_{n-2} + 1 + (k-2)(bn-b+1)\\ &= b\cdot S_{n-2} + 1 + (k-1)(bn-b+1).
  \end{align*}
  \item[(\ref{it:Bbipartite})] We proceed by induction on $m$. Proposition~\ref{prop:sRamsey}(\ref{it:star}) yields $r_{1}(K_{1,n},b)\leq \binom{b+1}{2}(n-1)+b+1 = \binom{b+1}{2}^1(n-1)+b\cdot S_0+1$, as desired.
  So assume that we know the statement for $m$ and we want to prove it for $m+1$.
  Since $\binom{b+1}{2}^{m+1}(n-1) + b\cdot S_m+1 = \binom{b+1}{2}\left(\binom{b+1}{2}^m(n-1)+b\cdot S_{m-1}+1-1\right)+b+1$,
  Proposition~\ref{prop:sRamsey}(\ref{it:gcone})
  guarantees that any $b$-iatlon graph with this number of vertices contains a constrained copy of $K_1+K_{m,n}$ with $\leq (m+1)$ distinct edge labels.
  This finishes the induction step, as $K_{m+1,n}$ is a subgraph of $K_1+K_{m,n}$ and we can apply Proposition~\ref{prop:sRamsey}(\ref{it:subgraph}).
  \item[(\ref{it:Bkbipartite})] 
  We again proceed by induction on $k$, the case $k=1$ being given by the previous point.
  In the considered situation,
  $\binom{b+1}{2}^m(n-1)+b\cdot S_{m-1}+1 + (k-1)(bm+n) \geq \binom{b+1}{2}^m(n-1)+b\cdot S_{m-1}+1 \geq r_m(K_{m,n},b)$.
  Therefore, Proposition~\ref{prop:sRamsey}(\ref{it:bdunion}) and the induction assumption yield the desired inequality
  \begin{align*}
   r_{km}(k\cdot K_{m,n},b)&\leq m+n+(b-1)m + r_{(k-1)m}((k-1)\cdot K_{m,n},b)\\ &\leq bm+n + \binom{b+1}{2}^m(n-1)+b\cdot S_{m-1}+1 + (k-2)(bm+n)\\
   &=\binom{b+1}{2}^m(n-1)+b\cdot S_{m-1}+1 + (k-1)(bm+n).\qedhere
  \end{align*}
 \end{description}
\end{proof}

\subsection{Optimality}\label{sec:optimality}
In this section we provide some lower bounds for the function $p(G,b)$. Because of the
prominent role of the stars in the proof of Proposition~\ref{prop:sRamsey} and Theorem~\ref{thm:bounds},
we shall focus on establishing bounds for $p_(K_{1,n},b)$. 
 For $b\leq 2$, we provide an optimal bound. 
\begin{lemma}\label{lem:lowerbound}
 For each $b\geq 1$ and $n\geq 0$ there 
 is a $b$-iatlon graph on $bn$ vertices
 that does not contain any constrained copy of $K_{1,n}$.
\end{lemma}
\begin{proof}
 Let $H$ be a disjoint union of $b$ complete graphs $K_n$. For each edge $e$ in this union we set $L_e:=2^{V(H)\setminus e}$. Then for any set $T$ of $b+1$ vertices,
 two of these vertices lie in the same copy of $K_n$ and hence are connected via an edge $e$, which has $T\setminus e$ as one of its labels. It follows that $H$ is a $b$-iatlon graph. Clearly $H$ does not contain any constrained copy of $K_{1,n}$, since any such copy has $n+1$ vertices.
\end{proof}
In some cases, the provided bound is optimal:
If $n=0$, then $bn+1=1$ and one vertex suffices to find a constrained copy of $K_{1,0}=K_1$.
If $n=1$, then by the definition of $b$-iatlon graph, $bn+1=b+1$ vertices suffice to find a constrained copy of $K_2$.
If $b=1$, then $bn+1=n+1$ vertices suffices to obtain a constrained copy of $K_{1,n}$, see also Proposition~\ref{prop:sRamsey} (\ref{it:noholes}).
We are now going to show that the same is also true for $b=2$.

\begin{lemma}\label{lem:improved}
 Any $2$-iatlon graph $H$ on $2n+1$ vertices 
 contains a constrained copy of $K_{1,n}$.
\end{lemma}
Computer experiments revealed that in many further cases, $bn+1$ vertices in $H$ suffice to force a constrained copy of $K_{1,n}$. 
We hence conjecture that $bn+1$ vertices always suffice, see~Conjecture~\ref{conj:one}.

Nevertheless, we do not know whether it holds in general.
\begin{proof}
 The lemma obviously holds for $n\in\{0,1\}$.
 Without loss of generality we may assume
 that all labels of $H$ have cardinality $1$.
 During the proof we write $xy$ to denote the edge $\{x,y\}$.
 Let $x,y,p\in V(H)$ be a triple of distinct points. If $\{p\}\in L_{xy}$,
 we say that the triple $xyp$ chooses $xy$.
 As $H$ is $2$-iatlon, each triple chooses at least one pair of its vertices. Without loss of generality, we assume that it chooses exactly one pair.
 
 Let $m$ be a maximal number of edges with the same label that are adjacent to a single vertex. If $m\geq n$,
 then $H$ contains a copy $C_0$ of $K_{1,n}$ where all edges have the same label $l$. In such a case $l$ has to be disjoint from $V(C_0)$ and hence $C_0$ is a constrained copy of $K_{1,n}$. 
 
 Therefore, we may assume that
 \begin{equation}
  m<n.\label{eq:smaller}
 \end{equation}
 
 We will now show that each pair $\{x,y\}$ of vertices is an edge of $H$ and has at least one label. So let $x,y\in V(H)$ be an arbitrary pair of vertices. If we look at all the triples $xyp$ for $p\in V(H)\setminus \{x,y\}$, at most $m$ of them choose $xp$ and at most $m$ of them choose $yp$.
 Therefore, the number of labels of $\{x,y\}$ is at least  $2n-1-2m\overset{\eqref{eq:smaller}}{\geq} 1$, and so $\{x,y\}$ is an edge. 
 By \eqref{eq:smaller} and the equation \eqref{eq:degreeBound} from the proof of Proposition~\ref{prop:sRamsey}(\ref{it:star}) (with $N=2n+1$ and $b=2$) also
 \begin{equation}
  \left\lceil\frac{2n-1}{3}\right\rceil\stackrel{\eqref{eq:degreeBound}}{\leq} m\stackrel{\eqref{eq:smaller}}{<} n.\label{eq:big}
 \end{equation}
 
 By our choice of $m$, there are two vertices $x$ and $y$ and $m$ edges $xx_1,xx_2,\ldots, xx_m$ of label $\{y\}$.
 It follows that these edges form a constrained copy $D$ of $K_{1,m}$. 
 As we already know, each pair $\{x,y\}$
 is an edge of $H$ with at least one label.
 Hence we may choose a label $\{z\}$ of $xy$. Note that $z$ is distinct from all $x_i$. Otherwise, it would mean that the triple $xyx_i$ choose both $y$ and $x_i$, contradicting our assumption. Set $S:=\{x,y, x_1,x_2,\ldots, x_m\}$, $T:=V(H)\setminus S$ and $R:=T\setminus \{z\}$,
 see~Figure~\ref{fig:onehole}.
 
 \begin{figure} 
  \begin{center}
  \includegraphics[page=4]{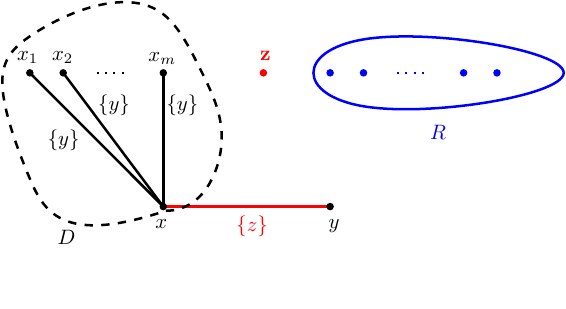}
  \caption{Proof for $b=2$}\label{fig:onehole}
  \end{center}
 \end{figure}

 We divide the rest of the proof into three cases.
 \begin{description}
  \item[1) $m=n-1$:] Here $|R|=n-1$.
  If any triple $xzr$ with $r\in R$ chooses $xz$, we may add $xz$ with label $\{r\}$ to $D$ and create a constrained copy of $K_{1,n}$; similarly if it chooses $xr$, we add $xr$ with label $\{z\}$. Hence we may assume that this is not the case and each such triple chooses $zr$ with label $\{x\}$.
  Consequently there are $n-1$ edges from $z$ with label $\{x\}$ and these form a constrained copy $D'$ of $K_{1,n-1}$.
  
  We now look at the triples $yx_iz$. If any of them chooses edge $zy$ or $zx_i$, this edge can be added to $D'$ to form a constrained copy of $K_{1,n}$. Consequently, we may assume that all these triples choose $yx_i$ with label $\{z\}$. Together with $xy$, we get $n$ edges with label $\{z\}$ incident to $y$,
  contradicting our choice of $m$. 
  Hence from now on, we may assume that
  \begin{equation}
   m \leq n-2.\label{eq:smallest}
  \end{equation}
  \item[2) There is $t\in T$ for which all labels of $xt$ lie in $S$.] 
  In this case we look at all the triples
  $xtt'$ with $t'\in T\setminus \{t\}$. By our assumptions none of them chooses $xt$ and at most $m$ choose $tt'$.
  Consequently at least $2n-2-m-m$ choose $xt'$ with label $\{t\}$. Since we may add all these edges to $D$ which already has $m$ edges, and since $m + 2n-2-m-m = 2n-2-m\overset{\eqref{eq:smallest}}{\geq} n$,
  we obtain a constrained copy of $K_{1,n}$.
  \item[3) Each edge $xt$, $t\in T$ has some label $l_{xt}\notin S$.] In this case we choose one such label for $xt$ and delete the remaining ones. Note that $l_{xt}\in T$.
  Then for each vertex $s\in T$ we set $S_s:=\{t\in T\mid e_t\text{ has label }s\}$.
  Since each edge $xt$ has a unique label, one of the labels, say $t_1$, is used at most once by the edges $xt$, $t\in T$. 
  We remove $t_1$, the vertex corresponding to the label of $xt_1$, and $S_{t_1}$ from $T$. In the remaining set at least one color, say $t_2$, is chosen at most once, etc.
  We continue until we remove all points.
  
  Since $|T|=2n-1-m \overset{\eqref{eq:big}}{\geq} 3(n-m-1)+1$, we obtain a set $Q$ of at least $(n-m)$ vertices $t_1,\ldots, t_{n-m}$ such that the labels of the edges $xt_1$, $xt_2$,\dots, $xt_{n-m}$ are in $T\setminus Q$.
  Consequently these edges with their labels can be added to $D$ to form a constrained copy of $K_{1,n}$.\qedhere
 \end{description} 
\end{proof}

\bigskip

\section{Proofs of the Radon and Helly theorems}\label{sec:radon}
 \subsection{Proof strategy}
 The proofs of Theorems~\ref{thm:radon}, \ref{thm:surfaces} and \ref{thm:manifolds} rely on the following fact, which we shall establish in the rest of the section.
 Suppose that we are given a finite simplicial complex $K$  and a sufficiently large point set $P$.
 Then it is possible to find a continuous map $f\colon |K|\to X$ and sets $P_\sigma\subseteq P$, $\sigma\in K$, such that for disjoint simplices $\sigma$ the sets $P_\sigma$ are disjoint and such that $f(|\sigma|)\subseteq \operatorname{cl} P_\sigma$.
 This fact shall be established via induction on $\dim K$.
 If $\dim K=0$ or $\dim K=1$, such maps are provided by~Proposition~\ref{prop:sRamsey}.
  In subsection~\ref{sec:twodimensions} we obtain~Theorem~\ref{thm:surfaces} as a corollary of the first induction step.
 In subsection~\ref{sec:final} we show how to finish the induction under homotopic assumptions.
 As a consequence we  obtain~Theorem~\ref{thm:manifolds}. The final step in the proof of~Theorem~\ref{thm:radon} consists in replacing the continuous maps in the proof by ``virtual maps'' (special chain maps between certain chain groups) and showing that even in this more general setting all the arguments go through. This is done in subsection~\ref{sec:virtual}.

\subsection{Toy example}\label{sec:examples}
Equipped with Theorem~\ref{thm:bounds}, we are now ready to prove an upper bound for our motivating example from the introduction.
\begin{proposition}\label{ex:drawing}
 Let $\operatorname{cl}$ be a closure operator on $\mathbb{R}^2$, such that for every $X\subseteq\mathbb{R}^2$ with $|X|=4$, $\operatorname{cl}(X)$ has at most $3$ path-connected components. Then
 $r(\operatorname{cl})\leq 562$.
\end{proposition}
\begin{proof}
 Given a set $S$ with $|S|\geq 562$,
 we choose distinct points $x_1,x_2,\ldots, x_{562}\in S$.
 We now consider a graph $H=(V,E)$,
 whose vertices are $x_1,\ldots, x_{562}$.
 Now given two distinct points $x_i,x_j$, it may happen that $x_i$ and $x_j$ lie in different components of $\operatorname{cl}\{x_i,x_j\}$.
 But if we take four distinct points $x_i,x_j,x_k,x_l$, then $\operatorname{cl}\{x_i,x_j,x_k,x_l\}$ has at most three path-connected components. Therefore, inside $\operatorname{cl}\{x_i,x_j,x_k,x_l\}$, two of the points can be connected with a path, see~Figure~\ref{fig:connection}.
 If for example, $x_i$ can be connected to $x_j$, we add the edge $\{x_i,x_j\}$ into $H$ and decorate it with the label $l_{ij}:=\{x_k,x_l\}$. (Note that one edge may be decorated with many different labels.)

 Since $H$ has $562=\binom{3+1}{2}^3\cdot 2 + 3\cdot (1 + \binom{3+1}{2}+\binom{3+1}{2}^2) + 1$ vertices and this value is larger or equal to $r_3(K_{3,3},3)$ by Theorem~\ref{thm:boundsColors}(\ref{it:Bbipartite}), in $H$ we can find a copy $C_0$ of $K_{3,3}$ and for each edge $e\in E(C_0)$ choose one of its labels $l_e$ such that
 \begin{enumerate}
  \item this label $l_e$ is disjoint from $V(C_0)$
  \item for disjoint edges $e,f\in E(C_0)$,
  the labels $l_e$ and $l_f$ are disjoint.
 \end{enumerate}
 Then we may draw $C_0$ in the plane in such a way that for every edge $e=\{x_i,x_j\}\in E(C_0)$ its image $x_{i,j}$ connects $x_i$ to $x_j$ and lies in $\operatorname{cl}(\{x_i,x_j\}\cup l_e)$. 
 Hanani-Tutte theorem~\cite{Chojnacki1934,Tutte70} then yields two disjoint intersecting edges. Without loss of generality, let $x_{1,2}$ intersect $x_{3,4}$. Then $x_{1,2}$ lies in $\operatorname{cl}\bigl(\{x_1,x_2\}\cup l_{12}\bigr)$ and $x_{3,4}$ lies in $\operatorname{cl}\bigl(\{x_3,x_4\}\cup l_{34}\bigr)$. So $S_1=\{x_1,x_2\}\cup l_{12}$ and $S_2=\{x_3,x_4\}\cup l_{34}$ are two disjoint sets with intersecting closures. Since $S$ was an arbitrary set with $|S|\geq 562$, this shows that $r(\operatorname{cl})\leq 562$, as desired.
\end{proof}

 \subsection{The two-dimensional case}\label{sec:twodimensions}
Using Theorem~\ref{thm:bounds}, we easily obtain bounds on Radon numbers for surfaces.
\begin{proof}[Proof of~Theorem~\ref{thm:surfaces}]
 Set $x:=\binom{b+1}{2}$ and
 \[N:=\binom{b+1}{2}^3(3-1)+b\left(\binom{b+1}{2}^2 + \binom{b+1}{2}+1\right)+1 + g(3b+3) = 
 2x^3 + bx^2 + bx + b + 3g(b+1).
 \] 
 Now, we proceed similarly as in the proof of Proposition~\ref{ex:drawing}.
 Let $S\subseteq M$ be an arbitrary set with $N$ or more points. We let $N$ of these points $x_1,x_2,\ldots, x_N$ to be vertices of $H$. 
 By the assumption, in each $(b+1)$-tuple $T$ of these points, some two of them can be connected with a path inside $\operatorname{cl} T$. If these vertices are $x$ and $y$, we put $\{x,y\}$ into $E(H)$ and decorate it with the label $l_{\{x,y\}}:=T\setminus\{x,y\}$.
 
 This way we obtain a $b$-iatlon graph on $N$ vertices.
 By~Theorem~\ref{thm:boundsColors}(\ref{it:Bkbipartite}), $N\geq r_{3(g+1)}((g+1)\cdot K_{3,3},b)$, 
 and so this $b$-iatlon graph contains a constrained copy $C_0$ of $(g+1)\cdot K_{3,3}$. 
 We draw each edge $e=\{x,y\}\in E(C_0)$ as a path that connects $x$ to $y$ inside $\operatorname{cl}\{e\cup l_0(e)\}$, where $l_0(e)$ denotes the chosen label of $e$ in the constrained copy $C_0$. The drawing is possible due to the definition of a constrained copy and our choice of the labels.
 
 By a result of Schaefer et al.~\cite{schaefer2013}, any drawing of $(g+1)\cdot K_{3,3}$ on $M$ contains two independent edges $e$ and $f$, whose images intersect. Since the image of $e$ lies in $\operatorname{cl}\{e\cup l_0(e)\}$ and the image of $f$ in $\operatorname{cl}\{f\cup l_0(f)\}$, this proves that $S_1:=e\cup l_0(e)$ and $S_2:=f\cup l_0(f)$ are two disjoint sets, whose closures intersect. Since $S$ was arbitrary, the theorem follows.
\end{proof}

\subsection{Higher dimensions}\label{sec:final}
To extend the results from graphs to higher
dimensions, that is, to simplicial complexes, we need an analogue of constrained graphs. This time we also include the information about the drawing of the simplicial complex.
The definition is a modification of the definition by Patáková~\cite{patakova2020} and Goaoc, Paták, Patáková, Tancer, Wagner~\cite{hb17}.
\begin{definition}[Constrained map]\label{d:constr_map}
 Let $K$ be a simplicial complex and
 let $\operatorname{cl}$ be some closure operator on a topological space $X$. Let $P\subseteq X$ be a set of points. Let $f\colon |K|\to X$ be a continuous  map. We say that $f$ is \emph{constrained by
  $(\operatorname{cl},\Phi)$} if $\Phi$ is a map from $K$ to $2^P$ that satisfies the following.
\begin{enumerate}[(i)]
 \item
 $\sigma\cap\tau=\emptyset\Rightarrow \Phi(\sigma)\cap\Phi(\tau)=\emptyset$.
 \label{it:(i)}
 \item $\sigma\subseteq\tau\Rightarrow \Phi(\sigma)\subseteq\Phi(\tau)$.\label{it:(ii)}
\item For all $\sigma \in K$, the image $f(|\sigma|)$ is contained in $\operatorname{cl}{\Phi(\sigma)}$.\label{it:(iii)}
\item For every vertex $v\in V(K)$, $\Phi(v)$ consists of one point only.\label{it:(iv)}
\end{enumerate}
If there is some map $\Phi\colon K\to 2^P$ such that the continuous map $f\colon |K|\to X$ is constrained by $(\operatorname{cl},\Phi)$,
we say that $f$ is \emph{constrained by $(\operatorname{cl}, P$)}.
If $\operatorname{cl}$ is clear from the context, we just say that $f$ is constrained by $\Phi$ or $P$, respectively.
\end{definition}
The main difference when compared to the older definitions by Goaoc et al.~\cite{hb17} and Patáková~\cite{patakova2020} is that we replace the conditions $\Phi(\emptyset)=\emptyset$ 
and $\Phi(\sigma\cap\tau)=\Phi(\sigma)\cap\Phi(\tau)$ with $\sigma\cap\tau=\emptyset\Rightarrow \Phi(\sigma)\cap\Phi(\tau)=\emptyset$.
This might seem like a small change, but
this change is crucial and underlies all the improvements presented in this paper. It namely allows to assign non-disjoint labels to non-disjoint edges, a fact that was heavily exploited in Proposition~\ref{prop:sRamsey}.

\begin{lemma}\label{lem:collision}
 Let $X$ be a topological space and $\operatorname{cl}$ be some closure operator on $X$.
 Let $K$ be a simplicial complex such that
 for every continuous map $f\colon K\to X$,
 there are two disjoint faces $\tau,\sigma\in K$ whose images intersect.
 If for every set $P$ of $N$ points,
 there is a continuous map $f\colon K\to X$ constrained by $P$, then $r(\operatorname{cl})\leq N$.
\end{lemma}
\begin{proof}
 To show that $r(\operatorname{cl})\leq N$,
 it suffices to show that for every set $P\subseteq X$ of $N$ points, there are two disjoint subsets of $P$, whose closures intersect. So let $P\subseteq X$ be an arbitrary set of $N$ points.
 Let $f\colon K\to X$ be a map constrained by $P$. By the assumption of the lemma,
 there are two disjoint faces $\sigma,\tau$,
 whose images intersect.
 Since $f$ is constrained by $(\operatorname{cl},\Phi)$, $f(|\sigma|)$ is contained in $\operatorname{cl}(\Phi(\sigma))$ and $f(|\tau|)$ is contained in $\operatorname{cl}(\Phi(\tau))$. Therefore, $\Phi(\sigma)$ and $\Phi(\tau)$ are two disjoint subsets of $P$, whose closures intersect. 
\end{proof}

\begin{lemma}\label{lem:induction}
 Let $X$ be a topological space and $\operatorname{cl}$ be some closure operator on $X$.
 Let $K$ be a $d$-dimensional simplicial complex $K$ on $n\geq 2$ vertices and $b$ be a positive integer. Set $S_{n-2}:=\sum_{j=0}^{n-2}\binom{b+1}{2}^j$. If
 \begin{itemize}
  \item for all $R$ of size $b+1$,
  $\operatorname{cl}(R)$ has at most $b$ path-connected components and
  \item for all $i\geq 1$, $i<d$ and all sets $R\subseteq X$ with $i+b+1\leq |R|\leq i+2 + (b-1)\binom{i+2}{2}$, $\pi_i(\operatorname{cl} R)=0$, 
  \item and if $P\subseteq X$ is a set with $|P|\geq b\cdot S_{n-2} + 1$,
  \end{itemize}
  then there exists a continuous map $f\colon |K|\to X$ constrained by $P$.
\end{lemma}
\begin{proof}
 Let the vertices of $K$ be $v_1,v_2,\ldots, v_n$. Let $|P|=N$. Note that $N\geq n$.
 We proceed by induction on $d$.
 Moreover, we construct $\Phi$
 in such a way that for each vertex $v$ of $K$, $\Phi(v)$ will consist of a single point, and for each face $\sigma$ with $\dim\sigma\geq 1$, we will have $|\sigma|+b-1\leq |\Phi(\sigma)|\leq |\sigma| + (b-1)\binom{|\sigma|}{2}$.
 
 If $d=0$, we choose $n$ pairwise distinct points $x_1,x_2,\ldots, x_n$, map $v_i$ to $x_i$ and set $\Phi(v_i)=\{x_i\}$.
 
 If $d=1$, $K$ is a graph. In this case, we consider a graph $H$ with $N$ distinct vertices $x_1,\ldots, x_N$ from $P$. By our assumptions, for each $(b+1)$-tuple $T$ of vertices of $H$, some two of them can be connected inside $\operatorname{cl}T$. If these vertices are $x$ and $y$, we add the edge $\{x,y\}$ to $H$ and decorate it with the label $l_{\{x,y\}}:=T\setminus\{x,y\}$. 
 This way we obtain a $b$-iatlon graph $H$ on $N$ vertices.
 Since $N\geq r_{n-1}(K_n,b)$ by~Theorem~\ref{thm:boundsColors}(\ref{it:five}), $H$ contains a constrained copy $C_0$ of $K_n$,
 and hence also a constrained copy $C_1$ of $K\subseteq K_n$. For each vertex $v\in V(C_1)$, we set $\Phi(v)=\{v\}$.
 For each edge $e\in E(C_1)$ we set $\Phi(e)=e\cup l_e$. Then $|\Phi(e)|=b+1$, as desired.
 Since $C_1$ is isomorphic to $K$ and since each edge $e=\{x,y\}\in E(C_1)$ comes with the induced path that connects $x$ to $y$ inside $\operatorname{cl}(e\cup l_e)$, we obtain a map $f\colon |K|\to X$ that is constrained by $\Phi$.

 If $d\geq 2$, then for each $k$-dimensional face $\sigma$ of $K$,  we set $\Phi(\sigma):=\{\Phi(v)\mid v\in\sigma\}\cup \bigcup_{\substack{e\subseteq \sigma\\ |e|=2}}l_e$, where $l_e$ comes from the case $\dim K=1$. In particular, $|\sigma|+b-1\leq |\Phi(\sigma)|\leq |\sigma|+(b-1)\binom{|\sigma|}{2}$, as desired. One easily checks that for disjoint faces $\sigma,\tau$ of $K$, $\Phi(\sigma)\cap \Phi(\tau)=\emptyset$ and that for $\tau\subseteq\sigma$, $\Phi(\tau)\subseteq\Phi(\sigma)$.
 
 We then proceed by induction.
 Hence assume that we already have a constrained drawing $f$ of the $(d-1)$-skeleton of $K$.  
 In particular, for each $d$-simplex $\sigma\in K$, we have a map $f\colon |\partial\sigma|\to\operatorname{cl}\Phi(\sigma)$.
 By our assumption $\pi_{d-1}(\operatorname{cl}\Phi(\sigma))=0$. But the equality
 $\pi_{d-1}(\operatorname{cl}\Phi(\sigma))=0$ just means that any map $f$ from the boundary $|\partial \sigma|$ of $\sigma$ to $\operatorname{cl}\Phi(\sigma)$ can be extended to a map from the whole simplex $|\sigma|$ to $\operatorname{cl}\Phi(\sigma)$.
 
 If we now perform this extension for all $d$-dimensional faces $\sigma$, we obtain a map $f\colon |K|\to M$ that is constrained by $\Phi$, which finishes the induction step.
\end{proof}

\begin{proof}[Proof of~Theorem~\ref{thm:manifolds}]
Let $L$ be a simplicial complex and $M$ a manifold.
A continous map $f\colon |L|\to M$ is called an \emph{almost-embedding}, if disjoint simplices of $L$ have disjoint images.
As was shown by Paták and Tancer~\cite[Theorem 1.2]{Patak_Tancer}, if $n$ and $k$ are arbitrary integers, and there is an almost embedding of the $k$-skeleton of an $n$-simplex into a compact PL manifold $N$, then $n\leq (2k+1)+(k+1)\beta_k(M;\mathbb{Z}_2)$.

 Let $K$ be the $d$-skeleton of the $\left(2d+2 + (d+1)\beta_d(M;\mathbb{Z}_2)\right)$-dimensional simplex. This complex has $t:=2d+3+(d+1)\beta_d(M;\mathbb{Z}_2)$ vertices.
 By Paták and Tancer's result mentioned above~\cite[Theorem 1.2]{Patak_Tancer}, for any continuous map $f\colon |K|\to M$, there are two disjoint simplices $\sigma,\tau\in K$ with $f(|\sigma|)\cap f(|\tau|)\neq \emptyset$.
 Set $S_{t-2}:=\sum_{j=0}^{t-2}\binom{b+1}{2}^j$.
 By Lemma~\ref{lem:induction}, for any set of $P\subseteq M$ with $N\geq b\cdot S_{t-2} + 1$, there is some continuous map $f\colon |K|\to M$ which is constrained by $(\operatorname{cl},P)$.
 Hence $r(\operatorname{cl})\leq N$ by~Lemma~\ref{lem:collision}.
\end{proof}

\subsection{Replacing homotopy with homology}\label{sec:virtual}
If we replace the condition $\pi_i(\operatorname{cl} R)=0$ with $\widetilde{H}_i(\operatorname{cl} R)=0$, then in the proof of Lemma~\ref{lem:induction}, we are only guaranteed that $f(|\partial\sigma|)$ can be filled by some singular chain. This way, instead of a proper map, we obtain a chain map $\varphi$ from the simplicial chains $\widetilde{C}_*^\Delta(K;\mathbb{Z}_2)$ to the singular chains $\widetilde{C}_*(X;\mathbb{Z}_2)$.
This chain map will map points to points.
In particular, $\varphi_{-1}$ will be an isomorphism. We call such chain maps \emph{non-trivial}.

Let us now explain the details.
A \emph{singular $k$-simplex} in a topological space $X$ is any continuous map 
$\gamma$ from the standard $k$-simplex $\operatorname{conv}\{e_1,\ldots, e_{k+1}\}$ into $X$. Its support $\operatorname{supp}(\gamma)$ is defined as the image of this map.
A \emph{singular chain} in $\widetilde{C}_*(X;\mathbb{Z}_2)$ is a formal linear combinations $\sum_{i=1}^n \gamma_i$ of singular simplices $\gamma_i$ with addition modulo $2$, i.e. $\gamma_i+\gamma_i=0$.
The support of such a chain is defined as $\bigcup_{i=1}^n \operatorname{supp}\gamma_i$.
A linear map $\varphi$ from 
$\widetilde{C}_*^\Delta(K;\mathbb{Z}_2)$ to $\widetilde{C}_*(X;\mathbb{Z}_2)$ is a chain, if it respects the boundary operator: 
$\varphi\circ \partial = \partial\circ \varphi$.

\begin{definition}[Constrained chain map]\label{d:constr_chain}
 Let $K$ be a simplicial complex and
 let $\operatorname{cl}$ be some closure operator on a topological space $X$. Let $P\subseteq X$ be a set of points. Let $\varphi\colon \widetilde{C}_*^\Delta(K;\mathbb{Z}_2)\to \widetilde{C}_*(X;\mathbb{Z}_2)$ be a non-trivial chain map. We say that $\varphi$ is \emph{constrained by
  $(\operatorname{cl},\Phi)$} if $\Phi$ is a map from $K$ to $2^P$ that satisfies the following.
\begin{enumerate}[(i)]
 \item
 $\sigma\cap\tau=\emptyset\Rightarrow \Phi(\sigma)\cap\Phi(\tau)=\emptyset$.
 \label{it:c(i)}
 \item $\sigma\subseteq\tau\Rightarrow \Phi(\sigma)\subseteq\Phi(\tau)$.\label{it:c(ii)}
\item For all $\sigma \in K$, $\operatorname{supp}\varphi(\sigma)\subseteq \operatorname{cl}{\Phi(\sigma)}$.\label{it:c(iii)}
\item For every vertex $v\in V(K)$, $\Phi(v)$ consists of one point only.\label{it:c(iv)}
\end{enumerate}
If there is some map $\Phi\colon K\to 2^P$ such that $\varphi$ is constrained by $(\operatorname{cl},\Phi)$, we say that $\varphi$ is \emph{constrained by $(\operatorname{cl}, P$)}.
\end{definition}

\begin{lemma}\label{lem:chaincollision}
 Let $X$ be a topological space and $\operatorname{cl}$ be some closure operator on $X$.
 Let $K$ be a simplicial complex such that
 for every non-trivial chain map $\varphi\colon \widetilde{C}_*^{\Delta}(K;\mathbb{Z}_2)\to \widetilde{C}_*(X;\mathbb{Z}_2)$, there are two disjoint faces $\tau,\sigma\in K$ with $\operatorname{supp} \varphi(\tau)\cap\operatorname{supp}\varphi(\sigma)\neq\emptyset$.
 If for every set $P$ of $N$ points,
 there is a non-trivial chain map $\varphi\colon \widetilde{C}_*^{\Delta}(K;\mathbb{Z}_2)\to \widetilde{C}_*(X;\mathbb{Z}_2)$ constrained by $P$, then $r(\operatorname{cl})\leq N$.
\end{lemma}
The proof is completely analogous to the proof of Lemma~\ref{lem:collision}.

\begin{lemma}\label{lem:chaininduction}
 Let $X$ be a topological space and $\operatorname{cl}$ be some closure operator on $X$.
 Let $K$ be a $d$-dimensional simplicial complex $K$ on $n\geq 2$ vertices and $b$ be a positive integer. Set $S_{n-2}:=\sum_{j=0}^{n-1}\binom{b+1}{2}^j$. If
 \begin{enumerate}
  \item for all $R$ of size $b+1$,
  $\operatorname{cl}(R)$ has at most $b$ path-connected components and
  \item for all $i\geq 1$, $i<d$ and all sets $R\subseteq X$ with $i+b+1\leq |R|\leq i+2 + (b-1)\binom{i+2}{2}$, $\widetilde{H}_i(\operatorname{cl} R;\mathbb{Z}_2)=0$, \label{it:chainrestriction}
  \item and if $P\subseteq X$ is a set with $|P|\geq b\cdot S_{n-2} + 1$,
  \end{enumerate}
  then there exists a non-trivial chain map $\varphi\colon \widetilde{C}_*^\Delta(K;\mathbb{Z}_2)\to \widetilde{C}_*(X;\mathbb{Z}_2)$ constrained by $P$.
\end{lemma}
\begin{proof}
The proof is completely analogous to the proof of Lemma~\ref{lem:induction}, the only real difference is that we replace ordinary maps with non-trivial chain maps. Nevertheless, to convince the reader that that are no issues connected to this change, let us spell out the technical details.

Let $G$ be the $1$-skeleton of $K$.
Then $G$ has at most $n$ vertices and so, 
by Lemma~\ref{lem:induction}, there is a continuous map $f\colon G\to X$ which is constrained by some $\Phi$. Moreover, we may require that for every vertex $v$ of $G$, $|\Phi(v)|=1$ and for every edge $e=uw$ of $G$, $|\Phi(e)|=b+1$ and $\Phi(u)\cup\Phi(w)\subseteq \Phi(uw)$.

We let $\varphi=f_\sharp$ be the induced chain map from $\widetilde{C}_*^\Delta(G;\mathbb{Z}_2)$ to $\widetilde{C}_*(X;\mathbb{Z}_2)$. 
That is, if $\sigma$ is a vertex or an edge of $G$, let $\iota$ denote the standard homeomorphism from the $0$ or $1$ dimensional standard simplex onto $|\sigma|$. Then we set 
$\varphi(\sigma):=f\circ \iota$ and extend this assignment linearly to the whole chain group $\widetilde{C}_*^\Delta(G;\mathbb{Z}_2)$.
Then $\varphi\colon \widetilde{C}_*^\Delta(G;\mathbb{Z}_2)\to \widetilde{C}_*(X;\mathbb{Z}_2)$ is non-trivial chain map that is  constrained by $\Phi$.

For each $\sigma\in K$,
we set \[\Phi(\sigma):=\{\Phi(v)\mid v\in\sigma\}\cup \bigcup_{\substack{e\subseteq \sigma\\ |e|=2}}\Phi(e).\]
As one easily checks, if $\sigma,\tau\in K$ are  disjoint, then
$\Phi(\sigma)\cap\Phi(\tau)=\emptyset$. Also observe that for each face $\sigma\in K$,
\[b-1+|\sigma|\leq |\Phi(\sigma)|\leq |\sigma|+\binom{|\sigma|}{2}(b-1).\]

Now for $k=2,3,\ldots,\dim K$,
we inductively build a non-trivial chain map $\varphi\colon\widetilde{ C}_*^\Delta(K^{(k)};\mathbb{Z}_2)\to \widetilde{C}_*(X;\mathbb{Z}_2)$ that is constrained by $\Phi$.
Here $K^{(k)}$ denotes the $k$-dimensional skeleton of $K$.

So assume that we have already built a non-trivial chain map $\varphi\colon \widetilde{C}_*^\Delta(K^{(k)};\mathbb{Z}_2)\to \widetilde{C}_*(X;\mathbb{Z}_2)$ that is constrained by $\Phi$.
Then for every $(k+1)$-dimensional simplex $\sigma\in K$,
$\operatorname{supp}\varphi(\partial\sigma)\subseteq \operatorname{cl}(\Phi(\sigma))$.
Since $\varphi$ is a chain map, $\partial\varphi(\partial\sigma)=
\varphi(\partial^2\sigma)=0$.
Therefore $\varphi(\partial\sigma)$ is a $k$-cycle in $\widetilde{C}_*(\operatorname{cl}(\Phi(\sigma)); \mathbb{Z}_2)$. But since $|\sigma|=k+2$,
$\widetilde{H}_{k}(\operatorname{cl}(\Phi(\sigma)); \mathbb{Z}_2)=0$ by the second assumption  (\ref{it:chainrestriction}) of the lemma.
In other words, every $k$-cycle in $\widetilde{C}_*(\operatorname{cl}(\Phi(\sigma);\mathbb{Z}_2)$ is a boundary of some singular chain. In particular, there is some singular chain $\theta_\sigma$ with $\operatorname{supp}\theta_\sigma\subseteq  \operatorname{cl}(\Phi(\sigma))$
and $\partial\theta_\sigma=\varphi(\partial\sigma)$.
If we now define $\varphi(\sigma):=\theta_\sigma$,
for all $k$-dimensional simplexes $\sigma\in K$, then $\partial\circ \varphi = \varphi\circ\partial$  and so 
$\varphi\colon \widetilde{C}_*^\Delta(K^{(k+1)};\mathbb{Z}_2)\to \widetilde{C}_*(X;\mathbb{Z}_2)$ is a non-trivial chain map that is constrained by $\Phi$, as desired.\end{proof}

\begin{proof}[Proof of Theorem~\ref{thm:radon}]
Due a result by Goaoc, Paták, Patáková, Tancer, Wagner~\cite[Corollary 14]{hb17}, 
if $\varphi$ is a chain map from the $\left\lceil d/2\right\rceil$-skeleton of $(d+2)$-dimensional simplex into $\mathbb{R}^d$ that maps points to points,  there are two disjoint faces $\sigma,\tau$ of this skeleton, whose supports intersect.
Together with Lemma~\ref{lem:chaincollision} and Lemma~\ref{lem:chaininduction}, this  proves~Theorem~\ref{thm:radon}.
\end{proof}

\section{Open problems}\label{sec:open}
We finish by stating several open problems.
\begin{conjecture}
 Every $b$-iatlon graph $H$ with at least $bn+1$ vertices contains a constrained copy of $K_{1,n}$.\label{conj:one}
\end{conjecture}
According to the first part of Proposition~\ref{prop:sRamsey}(\ref{it:gcone}),
an affirmative solution of this conjecture would
show that $p(K_1+G,b)\leq bp(G,b)+1$. By induction, this would mean that 
any $b$-iatlon graph $H$ with $b^{n-1}+b^{n-2}+\ldots + 1$ vertices contains a constrained copy of $K_n$. Consequently, this would decrease the bounds in Theorem~\ref{thm:radon} to $r(\operatorname{cl})\leq
b^{d+2}+b^{d+1}+\ldots + 1=O(b^{d+2})$,
in Theorem~\ref{thm:surfaces} to $r(\operatorname{cl})\leq 3b^3 + b^2 + b + 1 + g\cdot (b+1)=O(b^3)$ and similarly for other results. 

The conjecture is trivially true for $b=1$ and all values of $n$; or for $n\in\{0,1\}$ and all values of $b$. Lemma~\ref{lem:improved} shows that it also holds for $b=2$ and all values of $n$. Apart from that several other small cases have been verified by computer.

Solving this conjecture would provide the first step towards establishing the exact values of $p_G(b)$ for various graphs $G$.
With respect to the results mentioned
in this paper, the values for $G\in\{K_n, k\cdot K_n$, $K_{m,n}$ and $k\cdot K_{m,n}\}$ are of particular
interest.

\begin{conjecture}
For each manifold $M$, there exists some simplicial complex $K$, such that for every non-trivial chain map $\varphi\colon \widetilde{C}_*^\Delta(K;\mathbb{Z}_2)\to \widetilde{C}_*(M;\mathbb{Z}_2)$ from the simplicial chains of $K$ to the singular chains on $M$ there are two disjoint faces $\sigma,\tau$ with $\operatorname{supp}\sigma\cap\operatorname{supp}\tau\neq\emptyset$. 
\end{conjecture}
Proving this conjecture would allow us to replace the $\pi_i(\operatorname{cl} X)=0$ for $i=1,2,\ldots, \lceil d/2\rceil-1$ in Theorem~ \ref{thm:manifolds} with the weaker condition $\widetilde{H}_i(\operatorname{cl} X;\mathbb{Z}_2)=0$ and it would also imply that the main theorems from the papers by Goaoc, Patáková, Paták, Tancer, Wagner~\cite{hb17} and Patáková~\cite{patakova2020} can be extended to manifolds.

\section*{Author contributions}
This is a work of a single author.

\section*{Acknowledgments}
The research was supported by the Czech Science Foundation grant no. 22-19073S.
I would also like to thank Zuzana
Patáková for very careful and detailed proof reading and many valuable comments.

\section*{Financial disclosure}
None reported.

\section*{Conflict of interest}

The authors declare no potential conflict of interests.

\bibliography{rb}
\bibliographystyle{alpha}


%

\end{document}